# NONHYPERBOLIC DEHN FILLINGS
# ON HYPERBOLIC 3-MANIFOLDS

## Mario Eudave-Muñoz and Ying-Qing Wu

Abstract. In this paper we will give three infinite families of examples of nonhyperbolic Dehn fillings on hyperbolic manifolds. A manifold in the first family admits two Dehn fillings of distance two apart, one of which is toroidal and annular, and the other is reducible and $\partial$-reducible. A manifold in the second family has boundary consisting of two tori, and admits two reducible Dehn fillings. A manifold in the third family admits a toroidal filling and a reducible filling with distance 3 apart. These examples establish the virtual bounds for distances between certain types of nonhyperbolic Dehn fillings.

## §1. Introduction

Given a slope $r$ on a torus boundary component $T_0$ of a 3-manifold $M$, the *Dehn filling* of $M$ along the slope $r$, denoted by $M(r)$, is the manifold obtained by gluing a solid torus $V$ to $M$ along $\partial V$ and $T_0$ so that $r$ bounds a meridian disk on $V$. A manifold is *simple* if it is irreducible, $\partial$-irreducible, atoroidal, and anannular. Thus a simple manifold is either hyperbolic, or a small Seifert fiber space, or it would be a counter example to the Geometrization Conjecture. In particular, if $M(r)$ has nonempty toroidal boundary, then it is simple if and only if it is hyperbolic [Th]. A Dehn filling $M(r)$ is of type $S$ (resp. $D$, $T$, $A$) if $M(r)$ contains an essential $S^2$ (resp. $D^2$, $T^2$, $A^2$), so it is reducible (resp. $\partial$-reducible, toroidal, annular). The bound $\Delta(X, Y)$ is the least nonnegative number $n$ such that if $M$ is a hyperbolic manifold which admits two Dehn fillings $M(r_1), M(r_2)$ of type $X, Y$, respectively, then $\Delta(r_1, r_2) \leq n$. The bounds $\Delta(X, Y)$ have been established, via the work of many people, for all the 10 possible choices of $(X, Y)$; see [GW2] for more details.

In some cases, the upper bound of $\Delta(X, Y)$ is reached only by a few manifolds. For example, it was shown in [GW1] that if $M(r_1)$ is annular and $M(r_2)$ is toroidal, then $\Delta(r_1, r_2) \leq 3$ unless $M$ is one of three special manifolds, for which $\Delta(r_1, r_2)$ is 4 or 5; moreover, there are infinitely many manifolds which admit two such Dehn fillings with $\Delta(r_1, r_2) = 3$. Thus $\Delta(A, T) = 5$, but the "virtual bound" to be defined below is 3. Similarly for $\Delta(T, T)$, see [Go]. The main results of this paper are the following.

1991 *Mathematics Subject Classification.* Primary 57N10.

Wu's research at MSRI was supported in part by NSF grant #DMS 9022140.

Typeset by $\mathcal{A}\mathcal{M}\mathcal{S}$-TEX

**Theorem 0.1.** *There are infinitely many hyperbolic manifolds $M$ which admit two nonhyperbolic Dehn fillings $M(r_1)$ and $M(r_2)$, such that $M(r_1)$ is toroidal and annular, $M(r_2)$ is reducible and $\partial$-reducible, and $\Delta(r_1, r_2) = 2$.*

**Theorem 0.2.** *There are infinitely many hyperbolic manifolds $M$ with two torus boundary components, each of which admits two reducible Dehn fillings $M(r_1), M(r_2)$, with $\Delta(r_1, r_2) = 1$.*

**Theorem 0.3.** *There are infinitely many hyperbolic manifolds $M$ which admit two nonhyperbolic Dehn fillings $M(r_1)$ and $M(r_2)$, such that $M(r_1)$ is reducible, $M(r_2)$ is toroidal, and $\Delta(r_1, r_2) = 3$.*

These theorems follow immediately from Theorems 2.6, 3.6 and 4.2 below. The first example satisfying the conditions in Theorem 0.1 was found by Hayashi and Motegi [HM], and the first example as in Theorem 0.2 was given by Gordon and Litherland [GLi].

Similar to $\Delta(X, Y)$, we define the *virtual bound* $\Delta_v(X, Y)$ of distances between type $X$ and type $Y$ Dehn fillings to be the maximal integer $n$ such that there are infinitely many hyperbolic manifolds $M$ which admit two Dehn fillings $M(r_1), M(r_2)$ of type $X, Y$ respectively, with $\Delta(r_1, r_2) = n$. If no such infinite family exist, define $\Delta_v(X, Y) = 0$. Thus $\Delta_v(X, Y) \leq \Delta(X, Y)$. The above theorems and some known results determine the virtual bounds of distances between certain types of nonhyperbolic Dehn fillings. The following is a table of $\Delta_v(X, Y)$.

Table 1.1. Virtual bound $\Delta_v(X, Y)$

As we can see, except for $\Delta_v(A, A)$, all the other $\Delta(X, Y)$ have been completely determined. In the table, $\Delta_v(T, T)$ is determined by Gordon [Go], $\Delta_v(T, A)$ by Gordon and Wu [GW1]. The upper bounds of the other entries in Table 1.1 are the same as that in [GW2], and the lower bounds of them are determined by Theorem 0.1 for $\Delta_v(D, T)$, $\Delta_v(D, A)$, and $\Delta_v(S, A)$; by Theorem 0.3 for $\Delta_v(S, T)$; by Gabai [Ga] and Berge [Be] for $\Delta_v(D, D)$; by Gordon and Wu [GW1] for $\Delta_v(A, A)$; and by Gordon and Litherland [GLi] for $\Delta_v(S, S)$. Theorem 0.2 gives a stronger result about type S-S fillings, manely the manifolds can be chosen to have an extra torus boundary components. Also, it provides infinitely many examples of two essential planar surfaces in 3-manifolds with distinct boundary slopes, and one of the surfaces has unbounded number of boundary components.

We would like to thank Cameron Gordon and John Luecke for some interesting discussion on this topic.

## §2. Toroidal/annular fillings and reducible/$\partial$-reducible fillings

In this section we prove Theorem 2.6, which shows that there are infinitely many hyperbolic manifolds which admit two Dehn fillings of distance two apart, one of which is toroidal and annular, and the other is reducible and $\partial$-reducible. Let $Y = S^2 \times I$. Consider the tangles $T_p$ in $Y$ as shown in Figure 2.1, where a rectangle labeled by an integer $n$ denotes a rational tangle of slope $1/n$; in other words, it contains two vertical strings with $n$ left hand half twists.

Figure 2.1

Let $T_p(r)$ be the tangle obtained by filling the inside sphere $S_0$ of $Y$ with a rational tangle of slope $r$. The tangles $T_p(r)$ are drawn in Figure 2.2(a)–(d) for $r = \infty, 0, -1, -1/2$, respectively. From the pictures we have the following lemma. We use $T(r, s)$ to denote a Montesinos tangle consisting of two rational tangles associated to the rational numbers $r$ and $s$ respectively. See [Wu2] or [Mo1, Co] for more details about Montesinos tangles and algebraic tangles.

**Lemma 2.1.** *(1) $T_p(\infty)$ is the connected sum of a trivial tangle and a Hopf link.*

*(2) $T_p(0)$ is the Montesinos tangle $T[\frac{1}{2p-1}, \frac{-1}{2p+1}]$.*

*(3) $T_p(-1)$ is the Montesinos tangle $T[\frac{1}{2p+1}, \frac{-1}{2p-1}]$.*

*(4) $T_p(\frac{-1}{2})$ is an algebraic tangle obtained by summing a Montesinos tangle $T[\frac{1}{2p}, \frac{-1}{2p}]$ with a rational tangle $T[\frac{1}{2}]$. It is not a Montesinos tangle.*  $\square$

Figure 2.2

Let $M_p$ be the double branched covering of $Y$ with branch set the tangle $T_p$. Then $M_p$ is a compact orientable 3-manifold with boundary consisting of two tori $T_0$ and $T_1$, where $T_0$ is the lift of the inside sphere $S_0$. The $\infty$ and 0 slopes on $S_0$ lifts to a meridian-longitude pair on $T_0$, with respect to which the Dehn filling manifold $M_p(r)$ is the double covering of the 3-ball branched along the tangle $T_p(r)$. See [Mo2] for more details. Denote by $Q(r,s)$ the double branched cover of a Montesinos tangle $T[\frac{1}{r}, \frac{1}{s}]$. Note that when $|r|, |s| > 1$, $Q(r,s)$ is a Seifert fiber space with orbifold $D(r,s)$, which by definition is a disk with two cone points of angle $2\pi/|r|$ and $2\pi/|s|$. Denote by $C(r,s)$ the cable space of type $(r,s)$, that is, the exterior of a knot $K$ in a solid torus $V$ which is parallel to a curve on $\partial V$ representing $rl + sm$ in $H_1(\partial V)$, where $(m, l)$ is a meridian-longitude pair of $\partial V$. The above facts and Lemma 2.1 lead to the following lemma.

**Lemma 2.2.** *Suppose $p \geq 2$. The manifolds $M_p$ have the following properties.*

*(1) $M_p(\infty)$ is the connected sum of a solid torus and the projective space $RP^3$;*

*(2) $M_p(0) = Q(2p-1, -2p-1)$;*

*(3) $M_p(-1) = Q(2p+1, -2p+1)$;*

*(4) $M_p(-1/2)$ is a non Seifert fibered graph manifold containing a unique essential torus $T$, cutting it into a cable space $C(2,1)$ and a Seifert fiber space $Q(2p, -2p)$.*

*Proof.* (1) follows from the fact that the double branched cover of the Hopf link is $RP^3$, and connected sum of links and tangles downstairs corresponds to connected sum of manifolds upstairs. (2) and (3) follow from the definition of $Q(r,s)$.

To prove (4), notice that the Conway sphere in $T_p(-1/2)$ cutting off the tangle $T(2p, -2p)$ lifts to an essential torus $T$ upstairs. Since $T_p(-1/2)$ is not a Montesinos tangle, $M_p(-1/2)$ is not a Seifert fiber space, so $T$ is the Jaco-Shalen-Johannson decomposing torus because each side of it is a Seifert fiber space. Since each of $C(2, 1)$ and $Q(2p, -2p)$ are atoroidal, $T$ is the unique essential torus in $M_p(-1/2)$. □

It is easy to see that if $|p| \leq 1$ then $M_p$ is non hyperbolic. In the following, we will assume $M = M_p$ and $p \geq 2$, and show that $M$ is hyperbolic. Since $M$ has toroidal boundary, by [Th] we need only show that $M$ is irreducible, $\partial$-irreducible, non Seifert fibered, and atoroidal.

**Lemma 2.3.** *If $p \geq 2$, then $M$ is irreducible, $\partial$-irreducible, and non Seifert fibered.*

*Proof.* If $M$ is reducible, let $S$ be a reducing sphere. $S$ is separating, otherwise it would be a reducing sphere in all $M(r)$, contradicting Lemma 2.2(2). Let $W, W'$ be the two components of $M$ cut along $S$, with $W$ the one containing $T_0$. Let $\widehat{W'}$ be $W'$ with $S$ capped off by a 3-ball. Since $M(0)$ is the Seifert fiber space $Q(2p-1, -2p-1)$, which is irreducible, $W(0)$ must be a 3-ball, so $\widehat{W'} = M(0) = Q(2p-1, -2p-1)$. But then we have

$$M(\infty) = \widehat{W'} \# \widehat{W}(\infty) = Q(2p-1, -2p-1) \# \widehat{W}(\infty) \neq (S^1 \times D^2) \# RP^3,$$

which is a contradiction. Therefore $M$ is irreducible.

If $M$ is $\partial$-reducible, then after $\partial$-compression one of the $T_i$ becomes a sphere separating the two component of $\partial M$, hence is a reducing sphere, contradicting the above conclusion.

If $M$ is Seifert fibered, then $M(r)$ is Seifert fibered for all but at most one $r$, for which $M(r)$ is reducible. Since $M(-1/2)$ is irreducible and is not a Seifert fiber space, this is not possible. □

**Lemma 2.4.** *Suppose $T$ is an essential separating torus in an irreducible 3-manifold $M$, and suppose it is compressible in $M(r_1), M(r_2)$ with $\Delta(r_1, r_2) \geq 2$, where $r_i$ are slopes on $T_0 \subset \partial M$. Then $T$ and $T_0$ cobound a cable space in $M$, with cabling slope $r_0$ satisfying $\Delta(r_0, r_i) = 1$, $i = 1, 2$.*

*Proof.* Cut $M$ along $T$ and let $X$ be the component containing $T_0$. Then $T$ is compressible in $X(r_i)$ and $\Delta(r_1, r_2) \geq 2$, so by [Wu, Theorem 1] there is an essential annulus $A$ in $X$ with one boundary on $T$ and the other on $T_0$, with slope $r_0$, say. Since $T$ is essential in $M$, it is not parallel to $T_0$, so by [CGLS, Theorem 2.4.3] $T$ is compressible in $X(r)$ only if $\Delta(r_0, r) \leq 1$. We must have $\Delta(r_0, r_i) = 1$, because if $r_0 = r_1$ then we would have $\Delta(r_0, r_2) = \Delta(r_1, r_2) = 2$, a contradiction. Now the manifold $X(r_i)$ is homeomorphic to the manifold $Y$ obtained by cutting $X$ along $A$, so the torus component of $\partial Y$ corresponding to $T$ under the homeomorphism is compressible in $Y$. Since $M$ is irreducible, this implies that $Y$ is a solid torus. It follows that $X$ is a cable space with cabling slope $r_0$. □

**Lemma 2.5.** *$M$ is atoroidal.*

*Proof.* Assuming the contrary, let $T$ be an essential torus in $M$. Then $T$ must be separating, otherwise $M(r)$ would contain a nonseparating torus or, if $T$ becomes compressible in $M(r)$, a nonseparating sphere, for all $r$, which contradicts Lemma 2.2(1).

Let $W, W'$ be the two components of $M$ cut along $T$, with $W$ the one containing $T_0$. Since $M$ contains no nonseparating essential torus, by the Haken finiteness theorem (cf. [Ja, Page 49]), we may choose $T$ to be outermost in the sense that $W'$ contains no essential torus.

CLAIM. *$T$ is compressible in $M(-1/2)$.*

Recall from Lemma 2.2(4) that $M(-1/2)$ has a unique essential torus $T'$. So if $T$ is incompressible in $M(-1/2)$ then either it is boundary parallel or it is isotopic to $T'$. The first case is impossible, otherwise $T'$ would be an essential torus in $W'$, contradicting the choice of $T$. Therefore $T$ must be isotopic to $T'$ in $M(-1/2)$. It follows that either $W' = C(2,1)$, or $W' = Q(2p, -2p)$.

Since $M(0)$ is atoroidal, either $T$ is boundary parallel in $M(0)$ or it is compressible in $M(0)$. In the first case we would have $Q(2p-1, -2p-1) = M(0) = W' = C(2,1)$ or $Q(2p, -2p)$, which is absurd. In the second case let $D$ be a compressing disk of $T$ in $W(0)$, and let $\widehat{W'}$ be the manifold obtained by capping off the sphere boundary component of $W' \cup N(D)$ with a 3-ball. Then $\widehat{W'}$ is a summand of $M(0) = Q(2p-1, -2p-1)$, so either $\widehat{W'} = Q(2p-1, -2p-1)$ or $\widehat{W'} = S^3$. However, this is impossible whether $W' = C(2,1)$ or $W' = Q(2p, -2p)$ because $\widehat{W'}$ is obtained from $W'$ by Dehn filling on $T$ along certain slope, and it is easily seen that when $p \geq 2$ none of the Dehn fillings on such $W'$ could produce $Q(2p-1, -2p-1)$ or $S^3$. This completes the proof of the claim.

Since $M(\infty)$ contains no incompressible torus, $T$ is compressible in $M(\infty)$. By the claim above, $T$ is also compressible in $M(-1/2)$. Since $\Delta(\infty, -1/2) = 2$, it follows from Lemma 2.4 that $W$ is a cable space $C(p,q)$ with cabling slope $r_0$ satisfying $\Delta(r_0, \infty) = \Delta(r_0, -1/2) = 1$. Solving these equalities, we have $r_0 = 0$ or $-1$. Now we have $W(r_0) = L(p,q) \# (S^1 \times D^2)$, so $M(r_0)$ should have a lens space summand. On the other hand, we have shown that $r_0 = 0$ or $-1$, and in either case by Lemma 2.2 $M(r_0)$ is a prime manifold with torus boundary. This contradiction completes the proof that $M$ is atoroidal. $\square$

**Theorem 2.6.** *The manifolds $M_p$, $p \geq 2$, are mutually distinct hyperbolic manifolds, each admitting two nonhyperbolic Dehn fillings $M(r_1)$ and $M(r_2)$, such that $M(r_1)$ is toroidal and annular, $M(r_2)$ is reducible and $\partial$-reducible, and $\Delta(r_1, r_2) = 2$.*

*Proof.* Consider the manifolds $M_p$ which is the double cover of $Y = S^2 \times I$ branched along the tangle $T_p$ in Figure 2.1. By Lemmas 2.3 and 2.5, $M_p$ are hyperbolic for all $p \geq 2$. By Lemma 2.2, $M_p(\infty)$ is reducible and $\partial$-reducible, and $M_p(-1/2)$ is the union of $C(2,1)$ and $Q(2p, -2p)$ along a torus, hence is toroidal and annular because there is an essential annulus in $C(2,1)$ with both boundary on the outside torus $T_1$. Since $\Delta(\infty, -1/2) = 2$, $M_p$ satisfy all the conditions of the theorem. It remains to show that $M_p$ and $M_q$ are non homeomorphic when $p, q \geq 2$ and $p \neq q$.

Let $T_0$ (resp. $T_0'$) be the torus of $\partial M_p$ (resp. $\partial M_q$) on which the Dehn fillings are performed. Let $(m, l)$ (resp. $(m', l')$) be the meridian-longitude pair on $T$ (resp. $T'$) chosen as in Lemma 2.2. Let $f : M_p \to M_q$ be a homeomorphism.

There is a homeomorphism of $Y$ interchanging the two sphere boundary components, and leaving $T_p$ invariant, which induces a self homeomorphism of $M_p$ interchanging the two boundary components. This can be seen by redrawing the tangle in Figure 2.1 as in Figure 2.3(a), where the sphere $S_0$ represents the inside sphere in Figure 2.1, and $S_1$ the outside sphere. After an isotopy the picture becomes that in Figure 2.3(b). (Note that the isotopy have changed the position of the endpoints of the tangle on the spheres, but that does not matter.) Now blow up the sphere $S_0$, we get the same picture as that in Figure 2.1, with $S_0$ and $S_1$ interchanged. Thus without loss of generality we may assume that $f$ maps $T_0$ to $T_0'$.

Figure 2.3

Since $M_p(\infty)$ is $\partial$-reducible, by [Sch] $M_p(r)$ is irreducible for all $r \neq \infty$. Hence the reducing slope $\infty$ is unique, so $f$ must send $m$ to $m'$. Assume $f(l) = l' + km'$. Because of uniqueness of Seifert fibration, neither of $M_p(0)$ or $M_p(-1)$ is homeomorphic to $M_q(0)$ or $M_q(-1)$ when $p, q \geq 2$ and $p \neq q$. Hence $k \neq 0, \pm 1$. Now $f$ sends the slope $-1/2$ to $(2k-1)/2$, so both $M_q(-1/2)$ and $M_q((2k-1)/2)$ are toroidal. We have $\Delta(-1/2, (2k-1)/2) = |4k| \geq 8$. On the other hand, by [Go], this happens only if $M_q$ is the figure 8 knot complement or the Whitehead sister link complement. Since $M_q$ have two boundary components, this is impossible. $\square$

## §3. MANIFOLDS ADMITTING TWO REDUCIBLE DEHN FILLINGS

In this section we will show that there are infinitely many hyperbolic manifolds with two torus boundary components, each admitting two reducible Dehn fillings. Consider the tangles $T_p$ in $Y = S^2 \times I$ as shown in Figure 3.1, where, as in Figure 2.1, a rectangle labeled by an integer $n$ denotes a rational tangle of slope $1/n$.

Figure 3.1

As in Section 2, we denote by $M_p$ the double branched cover of $Y$ branched along $T_p$, and by $T_p(r)$ the tangle obtained by filling the inside sphere $S_0$ with a rational tangle of slope $r$. Then the Dehn filling manifold $M_p(r)$ is the double cover of $Y$ branched along $T_p(r)$. The tangles $T_p(\infty)$ and $T_p(0)$ are drawn in Figure 3.2(a)–(b). We can see that $T_p(\infty)$ is the connected sum of $T(1/2, -1/2)$ and a Hopf link, while $T_p(0)$ is the connected sum of a Montesinos tangle $T(1/2p, -1/2p)$ and a Hopf link. Recall that $Q(r, s)$ denotes the Seifert fiber space which double branch covers the tangle $T(1/r, 1/s)$, and the double branched cover of a Hopf link is the projective space $RP^3$. Therefore we have the following lemma.

Figure 3.2

**Lemma 3.1.** *The manifolds $M_p$, $p \neq 0$, have the following properties.*
  *(1) $M_p(\infty) = Q(2, -2) \# RP^3$;*
  *(2) $M_p(0) = Q(2p, -2p) \# RP^3$.*   □

Thus each $M_p$ admits two reducible Dehn fillings. In below, we will assume $M = M_p$, and $p \geq 2$. We need to show that $M$ is hyperbolic. Let $T_0$ be the component of $\partial M$ on which the Dehn fillings are performed. Thus $T_0$ covers the inside sphere $S_0$ in Figure 3.1. Let $T_1$ be the component of $\partial M$ covering the outside sphere $S_1$.

**Lemma 3.2.** *$M$ is irreducible.*

*Proof.* Assuming the contrary, let $S$ be a reducing sphere of $M$. Clearly $S$ is separating, otherwise $M(0)$ would contain a nonseparating reducing sphere, contradicting Lemma 3.1. Let $W, W'$ be the components of $M$ cut along $S$, with $W$ the one containing $T_0$. Denote by $\widehat{W}$ the manifold $W$ with the sphere boundary capped off by a 3-ball. Similarly for $\widehat{W'}$. Then $\widehat{W'}$ is a summand of both $M(0)$ and $M(\infty)$, so by Lemma 3.1 we must have $\widehat{W'} = RP^3$. This also shows that the reducing sphere in $M$ is unique up to isotopy, because if $S$ and $S'$ bound different punctured $RP^3$, then tubing them together would give a sphere which does not bound a punctured $RP^3$.

Let $\rho$ be the involution of $M$ which induces the branched covering. Since the reducing sphere $S$ is unique up to isotopy, by the equivariant sphere theorem [MSY], it can be chosen to be invariant under the involution $\rho$, hence it double branch covers a sphere $S'$ in the manifold $Y$ downstairs, which must cut off a 3-ball $B$ because one side of $S$ is $W'$, which does not contain the preimage of $S_0$ or $S_1$. Extending the involution $\rho|_S$ trivially over a 3-ball $D$, we get a double branched cover $\widehat{W'} \to S^3 = B \cup D'$, with branch set $L$ the union of $T' = T_p \cap B$ and a trivial arc in the attached 3-ball $D'$, which is the image of $D$ under the branched covering map. Since $\widehat{W'} = RP^3 = L(2,1)$, the link $L$ is the 2-bridge link associated to the number $1/2$, which is the Hopf link. Therefore, $T' = T_p \cap B$ is a tangle in $B$ consisting of an unknotted arc and a trivial circle $C$ around it.

We want to shown that no such pair $(B, T')$ exists in $(Y, T_p)$. Assuming the contrary, then $(B, T')$ would remain the same after filling the sphere boundaries $S_0, S_1$ of $Y$ with any rational tangles. The tangle $T_p$ has two circle components $C_1, C_2$, where $C_1$ denotes the one on the left in Figure 3.1. The circle component $C$ of $T'$ must be one of the $C_i$. However, after filling both $S_i$ with 0-tangle, $C_2$ has linking number $p \geq 2$ with one of the components of the resulting link, while after filling $S_0$ with 1-tangle and $S_1$ with $\infty$-tangle the circle $C_1$ has linking number 2 with one of the components of the resulting link, either case contradicting the fact that $C$ bounds a disk in $B$ intersecting the resulting link only once. □

**Lemma 3.3.** *$M$ is $\partial$-irreducible, and is not a Seifert fiber space.*

*Proof.* Since $\partial M$ consists of two tori, $M$ being $\partial$-reducible would imply that it is reducible, which would contradict Lemma 3.2. If $M$ is Seifert fibered (with two torus boundary components), then $M(r)$ would be reducible for at most one $r$, which would contradict Lemma 3.1. □

**Lemma 3.4.** *Let $W$ be an irreducible and $\partial$-irreducible 3-manifold. If both $W(r_1)$ and $W(r_2)$ are reducible and $\partial$-reducible, then $r_1 = r_2$.*

*Proof.* Let $T_0$ be the Dehn filling component of $\partial W$. Assume $r_1 \neq r_2$. Since $W(r_1)$ is $\partial$-reducible and $W(r_2)$ is reducible, by Scharlemann's theorem [Sch, Theorem 6.1], $r_2$ is a cabling slope, so there is an essential annulus $A_2$ in $W$ with boundary two copies of $r_2$ of opposite orientations. Similarly, we have an essential annulus $A_1$ in $W$ with boundary consisting of two copies of $r_1$ of opposite orientations. Isotope $A_1$ to intersect $A_2$ essentially. Then $A_1 \cap A_2$ consists of essential arcs on $A_i$, running from one boundary component to the other. By the parity rule on [CGLS, Page 279], if an arc component of $A_1 \cap A_2$ connects two components of $\partial A_1$ which have opposite orientations on $T_0$, then it must connect two components of $\partial A_2$ with the same orientation on $T_0$. This is a contradiction because the two boundary components of each $A_i$ have opposite orientations on $T_0$.   $\square$

**Lemma 3.5.** *$M$ is atoroidal.*

*Proof.* Consider an essential torus $T$ in $M$. Clearly $T$ is separating, otherwise $M(0)$ would contain a nonseparating torus or sphere, which would contradict Lemma 3.1. Let $W, W'$ be the two components of $M$ cut along $T$, where $W$ contains $T_0$. Note that $T$ cannot be boundary parallel in $M(0)$ or $M(\infty)$, otherwise $W'$, and hence $M$, would be reducible, which would contradict Lemma 3.2. Hence $T$ is compressible in both $W(0)$ and $W(\infty)$ because by Lemma 3.1 they are atoroidal. After compression, $T$ becomes a sphere in $W(0)$ and $W(\infty)$, so if $W$ contained $T_1$, then both $W(0)$ and $W(\infty)$ would also be reducible, which is impossible by Lemma 3.3. Hence we conclude that any essential torus in $M$ must separate the two boundary components of $M$.

Let $\rho : M \to M$ be the involution which induces the branch covering, and let $X$ be the fixed point set of $\rho$. Then $X$ covers the tangle $T_p$ in the manifold $Y$ downstairs. Since $T_p$ contains four arcs running from $S_0$ to $S_1$, $X$ has four arcs running from $T_0$ to $T_1$, hence each essential torus $T$ intersect $X$ at least four times.

By the equivariant torus theorem [MS, Theorem 8.6], there is a set of essential torus $\mathcal{T}$ in $M$ such that $\rho(\mathcal{T}) = \mathcal{T}$. Let $T$ be a component of $\mathcal{T}$. Since $X$ intersects $T$ in at least four points, we must have $\rho(T) = T$. Calculating the Euler number of $T/\rho$, we see that $X$ cannot intersect $T$ in more than four points. Hence $T$ intersects $X$ exactly four times, and $S = T/\rho$ is a sphere in $Y$ which intersects each of the four arc components of $T_p$ exactly once, and is disjoint from the circle components of $T_p$. Since the two circle components of $T_p$ have linking number 1, they must lie on the same side of $S$.

Let $Y_1, Y_2$ be the two components of $Y$ cut along $S$, with $Y_1$ the one disjoint from the circle components of $T_p$. Let $W_1, W_2$ be the components of $M$ cut along $T$, with $W_i$ covering $Y_i$. Consider the tangle $T_p'$ consisting of the arc components of $T_p$. Let $M'$ be the double cover of $Y$ branched along $T_p'$, let $T'$ be the torus in $M'$ that covers $S$, and let $W_i'$ be the part of $M'$ that covers $Y_i$. It can be seen from Figure 3.1 that $T_p'$ is isotopic to four straight arcs running from $S_0$ to $S_1$; hence $M' = T^2 \times I$. Since $T'$ is a torus separating the two components of $\partial M'$, it is isotopic to a horizontal torus $T^2 \times x$, so each $W_i'$ is also homeomorphic to $T^2 \times I$. Now we have $T_p \cap Y_1 = T_p' \cap Y_1$, therefore $W_1$, as the double cover of $Y_1$ branched along $T_p \cap Y_1$, is the same as $W_1'$, hence is a product $T^2 \times I$. But then $T$ is boundary parallel, contradicting the assumption that $T$ is an essential torus in $M$.   $\square$

**Theorem 3.6.** *The manifolds $M_p$, $p \geq 2$, are distinct hyperbolic manifolds, each admitting two reducible Dehn fillings $M(r_1), M(r_2)$ with $\Delta(r_1, r_2) = 1$.*

*Proof.* We have shown in Lemmas 3.1–3.5 that $M_p$ are hyperbolic manifolds admitting two reducible Dehn fillings $M_p(0)$ and $M_p(\infty)$, so it remains to show that the manifolds are all different.

Suppose $f : M_p \to M_q$ is a homeomorphism, $p > q \geq 2$. As in the proof of Theorem 2.6, it is easy to see that there is a self homeomorphism of $M_p$ interchanging the two boundary components, hence we may assume that $f$ maps $T_0$ to $T_0'$, where $T_0'$ and $T_1'$ are the boundary tori of $M_q$, with $T_0'$ the one covering the inside sphere.

By [GLu1], $M_i$ admits at most three reducible Dehn fillings, with mutual distance 1. Since $M_p(0) = Q(2p, -2p) \# RP^3$ is homeomorphic to neither $M_q(0)$ nor $M_q(\infty)$, $f$ maps the slope 0 to another reducing slope of $M_q$, which must be $\pm 1$ because it has distance 1 from 0 and $\infty$. Thus the only reducible Dehn filling of $M_q$ homeomorphic to $M_p(\infty)$ is $M_q(\infty)$, so $f$ sends the $\infty$ slope on $T_0$ to $\infty$ on $T_0'$. Similarly, it sends the $\infty$ slope on $T_1$ to $\infty$ on $T_1'$. Denote by $M_p(r, s)$ the manifold obtained by $r$ filling on $T_0$ and $s$ filling on $T_1$. Then we have $M_p(0, \infty) = M_q(\pm 1, \infty)$.

The manifold $M_k(r, s)$ is a double cover of $T_k(r, s)$, which is obtained from $T_k$ by filling the inside sphere with a rational tangle of slope $r$ and the outside sphere with one of slope $s$. One can check that $T_p(0, \infty)$ is the split union of a Hopf link and a trivial knot, while $T_q(\pm 1, \infty)$ is the connected sum of a Hopf link and a 2-bridge link associated to the rational number $\pm\frac{1}{4}$. Thus $M_p(0, \infty) = S^1 \times S^2 \# RP^3$, and $M_q(\pm 1, \infty) = L(4, \pm 1) \# RP^3$. Since these two manifolds are not homeomorphic, this is a contradiction.  $\square$

## §4. REDUCIBLE AND TOROIDAL FILLINGS

In this section we show that there are infinitely many hyperbolic manifolds which admit a reducible filling and a toroidal filling of distance 3 apart. Consider the tangles $T_p$ $(p \geq 2)$ in $Y$, as shown in Figure 4.1(a), where $Y$ is the 3-ball obtained by deleting the interior of the 3-ball $B$ in the figure from $S^3$. As before, let $T(r)$ be the union of $(Y, T_p)$ with a rational tangle of slope $r$, and let $M_p(r)$ be the double branched cover of $S^3$ branched along $T_p(r)$.

**Lemma 4.1.** *The manifold $M_p$ admits the following Dehn fillings.*

*(1) $M_p(\infty)$ is a non Seifert fibered, irreducible, toroidal manifold;*

*(2) $M_p(0)$ is a lens space $L((p-1)(p+3)+1, p+3)$;*

*(3) $M_p(1)$ and $M_p(1/2)$ are small Seifert fibered manifolds, but not lens spaces;*

*(4) $M_p(1/3) = L(3, 1) \# L(2, 1)$.*

Figure 4.1

*Proof.* The tangles $T(\infty), T(0), T(1), T(1/2), T(1/3)$ are shown in Figure 4.1(b)–(f), respectively. We can see that $T(\infty)$ is the union of $T[\frac{1}{2}, \frac{1}{-(p+2)}]$ and $T[\frac{1}{2}, \frac{1}{p}]$, and is not a Montesinos link; $T(0)$ is a 2-bridge link associated to the rational number $1/((p-1) + 1/(p+3)) = (p+3)/((p+3)(p-1)+1)$; $T(1)$ and $T(1/2)$ are Montesinos links consisting of three rational tangles; and $T(1/3)$ is the connected sum of a trefoil knot and a Hopf link. The result now follows by taking the double cover of $S^3$ branched along the corresponding links.   $\square$

**Theorem 4.2.** *The manifolds $M = M_p$, $p \geq 2$, are mutually distinct hyperbolic manifolds, each admitting two Dehn fillings $M(r_1)$ and $M(r_2)$, such that $M(r_1)$ is reducible, $M(r_2)$ is toroidal, and $\Delta(r_1, r_2) = 3$.*

*Proof.* Let $r_1 = 1/3$, and $r_2 = \infty$. Then $\Delta(r_1, r_2) = 3$, and by Lemma 4.1, $M(r_1)$ is reducible, $M(r_2)$ is toroidal. We need to show that $M_p$ are hyperbolic and mutually distinct.

$M$ is irreducible, otherwise a closed summand would survive after all Dehn fillings; but since $M(0)$ and $M(1)$ are non homeomorphic prime manifolds, this is impossible. $M$ is not a Seifert fiber space because two Dehn fillings $M(\infty)$ and $M(1/3)$ are non Seifert fibered. These imply that $M$ is $\partial$-irreducible. To prove $M$ is hyperbolic, it remains to show that $M$ is atoroidal.

If $T$ is an essential annulus in $M$, then it is compressible in $M(0)$, $M(1)$, $M(1/2)$ and $M(1/3)$. Since $M(0)$ is irreducible, $T$ must be separating. Let $W, W'$ be the components of $M$ cut along $T$, with $W$ the one containing $T_0$. Since $\Delta(1, 1/3) = 2$, by Lemma 2.4, $W$ is a cable space $C(r, s)$, with cabling slope $r_0$ satisfying $\Delta(r_0, 1) = \Delta(r_0, 1/3) = 1$. Solving these equalities, we have $r_0 = 0$ or $1/2$; but since $M(r_0)$ contains a lens space $L(r, s)$, we must have $r_0 = 0$. Let $\delta_0$ and $\delta_1$ be the slopes on $T$ which bound disks in $W(0)$ and $W(1/3)$, respectively. Since $0$ is the cabling slope, we have $\Delta(\delta_0, \delta_1) = |r| > 1$. Now $W(0)$ is the connected sum of a solid torus and $L(r, s)$, while $W(1/3)$ is a solid torus, so we have

$$M(0) = L(r, s) \# W'(\delta_0),$$
$$M(1/3) = W'(\delta_1).$$

Comparing the first equation with Lemma 4.1(2), we see that $W'$ is the exterior of a knot in $S^3$ with $\delta_0$ the meridional slope. But then since $\Delta(\delta_0, \delta_1) > 1$, by [GLu2] the manifold $M(1/3)$ would be irreducible, which would contradict Lemma 4.1(4). This completes the proof that $M$ is atoroidal, and hence hyperbolic.

It remains to show that the manifolds $M_p$ are mutually distinct. Assume there is a homeomorphism $f : M_p \cong M_q$, $p > q \geq 2$. Let $(m, l)$ and $(m', l')$ be the meridian-longitude pair of $M_p$ and $M_q$, respectively. By [CGLS], [GLu1] and [BZ, Theorem 0.1], a hyperbolic manifold admits a total of at most three reducible or cyclic Dehn fillings, with mutual distance 1. Thus two of the four slopes $0, 1/3, f(0), f(1/3)$ on $\partial M_q$ must be the same. But since $M_p(0)$ is not homeomorphic to $M_q(0)$ or $M_q(1/3)$, we must have $f(1/3) = 1/3$, and $f(0)$ is of distance 1 from 0 and $1/3$, so $f(0) = 1/2$ or $1/4$. The first is impossible because $M_q(1/2)$ is not a lens space. Hence $f(0) = 1/4$. Now $f(m) = f((m + 3l) - 3l) = (m' + 3l') \pm 3(m' + 4l')$, and we have $\Delta(m', f(m)) \geq 9$. Since both $m'$ and $f(m)$ are toroidal Dehn fillings slopes on $\partial M_q$, this contradicts [Go]. $\square$

## References

[Be]    J. Berge, *The knots in $D^2 \times S^1$ with nontrivial Dehn surgery yielding $D^2 \times S^1$*, Topology Appl. **38** (1991), 1–19.

[BZ]    S. Boyer and X. Zhang, *The semi norm and Dehn filling*, Preprint.

[Co]    J. Conway, *An enumeration of knots and links, and some of their algebraic properties*, Computational problems in abstract algebra, New York and Oxford: Pergamon, 1970, pp. 329–358.

[CGLS]  M. Culler, C. Gordon, J. Luecke and P. Shalen, *Dehn surgery on knots*, Annals Math. **125** (1987), 237–300.

[Ga]    D. Gabai, *Surgery on knots in solid tori*, Topology **28** (1989), 1–6.

[Go]    C. Gordon, *Boundary slopes of punctured tori in 3-manifolds*, Trans. Amer. Math. Soc. (to appear).

[GLi]   C. Gordon and R. Litherland, *Incompressible planar surfaces in 3-manifolds*, Topology Appl. **18** (1984), 121–144.

[GLu1]  C. Gordon and J. Luecke, *Reducible manifolds and Dehn surgery*, Topology **35** (1996), 385–409.

[GLu2]  ———, *Only integral Dehn surgeries can yield reducible manifolds*, Math. Proc. Camb. Phil. Soc. **102** (1987), 97–101.

[GW1]   C. Gordon and Y-Q. Wu, *Toroidal and annular Dehn fillings*, MSRI Preprint #1997-056.

[GW2]   ———, *Annular and boundary reducing Dehn fillings*, (to appear).

[HM]    C. Hayashi and K. Motegi, *Dehn surgery on knots in solid tori creating essential annuli*, Trans. Amer. Math. Soc. (to appear).

[Ja]    W. Jaco, *Three manifold topology*, Regional Conference Series in Math., vol. 43, 1980.

[MS]    W. Meeks and P. Scott, *Finite group actions on 3-manifolds*, Invent. Math. **86** (1986), 287 –346.

[MSY]   W. Meeks, L. Simon and S-T. Yau, *Embedded minimal surfaces, exotic spheres, and manifolds with positive Ricci curvature*, Annals Math. **116** (1982), 621–659.

[Mo1]   J. Montesinos, *Una familia infinita de nudos representados no separables*, Revista Math. Hisp.-Amer **33** (1973), 32–35.

[Mo2]   ———, *Surgery on links and double branched covers of $S^3$*, Annals of Math. Studies **84**, 1975, pp. 227–260.

[Sch]   M. Scharlemann, *Producing reducible 3-manifolds by surgery on a knot*, Topology **29** (1990), 481–500.

[Th]    W. Thurston, *Three dimensional manifolds, Kleinian groups and hyperbolic geometry*, Bull. Amer. Math. Soc. **6** (1982), 357–381.

[Wu1]   Y-Q. Wu, *Incompressibility of surfaces in surgered 3-manifolds*, Topology **31** (1992), 271–279.

[Wu2]   ———, *The classification of nonsimple algebraic tangles*, Math. Ann. **304** (1996), 457–480.

Mario Eudave-Muñoz, Instituto de Matematicas, UNAM, circuito Exterior, Ciudad Universitaria, 04510 Mexico D.F., MEXICO
    *E-mail address*: `eudave@servidor.unam.mx`

Ying-Qing Wu Department of Mathematics, University of Iowa, Iowa City, IA 52242
    *E-mail address*: `wu@math.uiowa.edu`